\documentclass[10pt]{article}
\usepackage[utf8]{inputenc}
\usepackage[T1]{fontenc}
\usepackage{amsmath}
\usepackage{amsfonts}
\usepackage{amssymb}
\usepackage{stmaryrd}
\usepackage{amsthm}

\newtheorem{thm}{Theorem}

\title{ Proving Taylor's Theorem from the Fundamental Theorem of Calculus by Fixed-point Iteration }

\author{Chris Thron, Texas A\&M University - Central Texas, Killeen TX}
\date{}

\begin{document}
\maketitle

\newcommand{\DD}{\ensuremath{\mathcal{D}} }
\newcommand{\II}{\ensuremath{\mathcal{I}_a} }
\newcommand{\LL}{\ensuremath{\mathcal{L}} }
\newcommand{\JJ}{\ensuremath{\mathcal{J}} }
\newcommand{\KK}{\ensuremath{\mathcal{K}} }
\renewcommand{\AA}{\ensuremath{\mathtt{S}} }
\newcommand{\uu}{\ensuremath{\mathbf{u}} }
\newcommand{\vv}{\ensuremath{\mathbf{v}} }
\newcommand{\oo}{\ensuremath{\mathbf{1}} }

\newcommand{\RR}{\ensuremath{\mathbb{R}} }

\begin{abstract}
Taylor's theorem (and its variants) is widely used in several areas of mathematical analysis, including numerical analysis, functional analysis, and partial differential equations. 
This article explains how Taylor's theorem in its most general form  can be proved simply as an immediate consequence of the Fundamental Theorem of Calculus (FTOC).   
The proof shows the deep connection between the Taylor expansion and fixed-point iteration, which is a foundational concept in numerical and functional analysis. One elegant variant of the proof also demonstrates the use of combinatorics and symmetry in proofs in mathematical analysis. Since the proof emphasizes concepts and techniques that are widely used in current science and industry, it can be a valuable addition to the undergraduate mathematics curriculum.

\end{abstract}


\section{Introduction}

Taylor’s theorem is one of the most important results that are taught in basic calculus classes. Its importance is both theoretical and practical. Taylor's theorem is a foundational result in the field of numerical analysis: many error estimates for numerical solutions to algebraic or differential equations are based on the Taylor expansion of the solution. Taylor series is also fundamental to the theory and application of differential equations, in which series solutions play a large role. In fact, the Cauchy-Kowalesky theorem that establishes the existence  and uniqueness of solutions to partial differential equations relies on Taylor expansion of the solution \cite{folland1995}.

It has been noted  that the proofs of Taylor’s theorem given in many textbooks are not well motivated \cite{zvonimir1990}. Most often, the theorem is derived either using Rolle's theorem \cite{thomas2013}, or as a consequence of the mean value theorem \cite{leithold1996}, or by repeated integration by parts \cite{courant1999}. A less common proof uses induction \cite{brauer1987}.
It would seem that such a key theorem should have a stronger motivation, and should not be derived as an incidental result from other theorems. 

The purpose of this paper is to show that Taylor’s theorem is in fact an immediate consequence of the fundamental theorem of calculus; and furthermore, the proof is a straightforward application of fixed point iteration, which is one of the most basic techniques in numerical and functional analysis. We thus both provide motivation for Taylor's theorem, and show its deep relationship with other areas of mathematical analysis. 

The paper is organized as follows. Section \ref{Background} gives necessary background for the mathematical concepts and notation used in the proof. Section \ref{Method} presents the new proof of Taylor's theorem. Finally, Section \ref{Discussion} discusses the possible instructional uses for this new proof.

\section{Background}\label{Background}
This section introduces notation and concepts necessary for understanding the proof of Taylor's theorem in Section~\ref{Method}. These concepts are indispensible in modern applied mathematics.
\subsection{Theorem statements}
The fundamental theorem of calculus expresses the inverse relationship between differentiation and integration. There are many alternative statements of the theorem that differ slightly from one another.  We will use the following formulation \cite{barcenas2000}:

\begin{thm}\label{thm0}{(Fundamental theorem of calculus for the Lebesgue integral)~~} 	
A function $f : [a, b] \
\rightarrow \RR$ is absolutely continuous if and only if it is
differentiable almost everywhere, its derivative $f' \in L^1[a,b]$ and , for each
$x \in [a, b]$,
\begin{equation}\label{FTOC}
f(x)=f(a)+\int_{a}^{x} f^{\prime}\left(t \right) d t .
\end{equation}
\end{thm}
Several alternative proofs of this theorem are cited in \cite{barcenas2000}.

Taylor's theorem may be stated as follows:
\begin{thm}\label{thm1}{(Taylor's theorem)~~} 
Given a  function $f : [a, b] \rightarrow \RR$ such that the $N$th derivative $f^{(N)}$ is absolutely continuous on $[a,b]$ , then for all $x \in [a,b]$ 
\begin{equation}\label{taylor}
\begin{aligned}
f(x)&=f(a)+(x-a)f'(a) + \frac{(x-a)^2}{2!}f''(a) + \ldots \\
& \qquad \ldots+ \frac{(x-a)^N}{N!}f^{(N)} + R_n(x),
\end{aligned}
\end{equation}
where $R_n(x)$ satisfies the inequality
\begin{equation}\label{boundR}
|R_n(x)| \le \sup_{c \in [a,x]} |f^{(N+1)}(c)| \frac{(x-a)^{N+1}}{(N+1)!},
\end{equation}
and $R_n(x)$ is given exactly by the expression
\begin{equation}\label{exactR}
R_n(x) = \int_a^x \frac{(x-t)^{N}}{N!}f^{(N+1)}(t) dt. 
\end{equation}

\end{thm}

\subsection{Linear operators and some basic properties}\label{ops}
The concept of a \emph{linear operator} is  not typically mentioned in basic calculus classes: however, there are good reasons for changing this practice. Modern computational mathematics is fundamentally based on the connection between linear algebra and analysis, and the idea of linear operator is a key aspect of this connection.
Linear operators play the same role in calculus that matrices play in linear algebra, as we shall now explain in more detail.

Let $V,W$ be finite-dimensional real vector spaces, and let $L:V\rightarrow W$ be a mapping (i.e. a function) with domain $V$ and codomain $W$.  The mapping $L$ is \emph{linear} if for all  vectors $\uu,\vv \in  V$ and constants $a,b \in \mathbb{R}$ we have  
\begin{equation}\label{linearDef}
L(a\uu + b\vv) = aL(\uu) + bL(\vv) 
\end{equation}
It is a well-known result in linear algebra that every linear mapping $L$ between finite-dimensional vector spaces can be associated with a unique matrix (which we denote here as $M$)  such that $L(\vv) = M\vv$. 

In order to draw the connection between matrices and linear operators, we observe that function spaces (such as $L^1[a,b]$ or the set of absolutely continuous functions) are infinite-dimensional vector spaces.  Therefore given two function spaces $A,B$, we may consider the set of mappings from $A$ to $B$. Such mappings are called \emph{operators}. An operator $\LL:A \rightarrow B$ is \emph{linear} if for any functions $f,g \in A$ we have  
\begin{equation}\label{linearOp}
\LL(af + bg) = a\LL(f) + b\LL(g), 
\end{equation}
which obviously corresponds exactly with \eqref{linearDef}. Because of the close connection between matrices and linear operators, the action of a linear operator is often written as if it were a multiplication: for example, $\LL(f)$ is written instead as $\LL f$, so that \eqref{linearOp} is written as 
\begin{equation}\label{linearOp2}
\LL(af + bg) = a\LL f + b\LL g. 
\end{equation}
Note that $\LL f$ is in itself a real-valued function, which is a mapping from \RR to \RR: thus $\LL f(x)$ denotes the value of the function $\LL f$ evaluated at the point $x$. 

Two important examples of operators are the  \emph{differentiation } (or \emph{derivative}) \emph{operator}
\begin{equation}\label{Ddef}
\DD f(x) := \lim_{\delta \rightarrow 0} \frac{f(x+\delta)-f(x)}{\delta}.
\end{equation}
where  \DD is defined on the set of absolutely continuous functions; and the \emph{integral operator} 
\begin{equation}\label{Idef}
\II g(x) :=  \int_a^x g(t) dt,
\end{equation}
which is defined for all $g \in L^1[a,b]$.

  Another operator that is defined for absolutely continuous functions is the \emph{evaluation operator}:
\begin{equation}
f \rightarrow f(a)\oo,
\end{equation}
where $a$ is any number in the domain of $f$, and \oo denotes the constant function that takes the value 1 for all values of $x$ in the domain of $f$. 

Operators can be composed to produce other operators. In the subsequent discussion, the \emph{associativity} of operator composition will play an important role: if $\JJ, \KK, \LL$ are compatible operators, then for any function $f$ in the domain of \LL we have
\begin{equation}
\JJ(\KK \LL f) = (\JJ\KK) \LL f.
\end{equation} 
It follows that parentheses  are unnecessary when writing operator compositions: for example, we may write $\II (\II (\DD (\DD f)))$ as $\II^2\DD^2f$. 

A linear combination of linear operators is also a linear operator. In Section~\ref{Method} we will make use of the following linear operator, which is defined as a linear combination:
\begin{equation}\label{ldef}
\LL f := f(a)\oo + \II \DD f,
\end{equation}
$\LL f$ is defined on the set of absolutely continuous real-valued  functions defined on any interval containing $a$.

 In Section~\ref{Method} we will make use of the \emph{monotonicity} property of the integral operator: if $g, h \in  L^1[a,x]$, then
\begin{equation}\label{Imon}
g \le h \implies \II g  \le \II h,
\end{equation}
where the inequalities in \eqref{Imon} hold pointwise throughout the interval $[a,x]$.  Since  
\begin{equation}\label{Imon0}
\begin{aligned}
- \sup_{c\in [a,x]}|g(c)| \cdot \oo \le g  \le \sup_{c\in [a,x]}|g(c)| \cdot \oo\\ 
\end{aligned}
\end{equation}
it follows from the linearity of \II that
\begin{equation}\label{Imon1}
\begin{aligned}
- \sup_{c\in [a,x]}|g(c)| \cdot \II \oo &\le \II g  \le \sup_{c\in [a,x]}|g(c)| \cdot \II \oo \\
  \implies  | \II g | &\le \sup_{c\in [a,x]}|g(c)| \cdot \II  \oo.
\end{aligned}
\end{equation}
By iterating \eqref{Imon}  $n$ times and applying $\II^n$ to \eqref{Imon0}, we may generalize \eqref{Imon1} to
\begin{equation}\label{Imon2}
\begin{aligned}
| \II^n g | \le \sup_{c\in [a,x]}|g(c)| \cdot \II^n  \oo.
\end{aligned}
\end{equation}

\subsection{Fixed points and fixed-point iteration}\label{Fixed}

A \emph{fixed point} is any mathematical object $x$ that satisfies an equation of the form $x = g(x)$, where $g$ is a function.  The concept of fixed point play a hugely important role in mathematical analysis, both theoretically and computationally. The online resource Mathworld  lists nine prominent fixed-point theorems that appear in diverse areas of mathematics, including the famous Brouwer and Banach fixed point theorems for topological and metric spaces respectively \cite{enwiki:1087275378}.  

Fixed point ideas also figure in other important theoretical results.   
The central limit theorem in probability is connected to the fact that the standard normal density function is a fixed point of the modified convolution operator $f(x) \rightarrow f(\sqrt{2}x)*f(\sqrt{2}x)$ \cite{feller1968}. Given a time-dependent differential equation $\frac{df}{dt} = A(f,x)$, then a steady-state solution satisfies $A(f,x)=0$ which can also be written as a fixed point equation: $f = f + A(f,x)$. Furthermore, periodic solutions correspond to fixed points of the Poincare map \cite{angenent2016}.  In general, any algebraic equation can be written as a fixed point equation: for example, we have for any real-valued function $f$
\begin{equation}
f(x) = 0 \iff x = x + f(x).
\end{equation}

On the practical side,  one of the most versatile and widely-applied techniques in numerical analysis is \emph{fixed point iteration}.
One example is Newton’s method for finding roots of functions. This method is  based on the fact that any root $x^*$ of a differentiable function $f:\mathbb{R}\rightarrow \mathbb{R}$,   
is a locally stable fixed point of the function $g(x) := x - \frac{f(x)}{f^{\prime}(x)}$ as long as $f^{\prime}(x^*) \neq 0$. 
As a result, $x^*$ can be estimated numerically by iteration: that is, given an initial guess $x_0$ that is sufficiently close to $x^*$ then the series $g(x_0), g(g(x_0)), g(g(g(x_0))), \ldots g^{(n)}(x_0) \ldots$ converges to $x^*$. 
Another example is the power method for finding the largest eigenvector-eigenvalue pair of a matrix (or more generally, a linear operator).  
Given a matrix $M$, the unit eigenvector \vv corresponding to the eigenvalue with largest absolute value is a solution to the fixed-point equation   $\vv = F(\vv)$ where $F(\vv) := \frac{M\vv}{|M\vv|}$. It follows that \vv can be estimated numerically by repeatedly iterating the function $F$:  $F^{(n)}(\vv) \xrightarrow[n \rightarrow \infty]{} \vv^*$ \cite{panju2011}.  In physics, perturbation series are often obtained through fixed-point iteration. A prominent example is the \emph{Born series} in electromagnetic and quantum scattering \cite{newton2013}.  

\section{Method}\label{Method}
\subsection{Basic proof of Taylor's theorem}
It is not commonly recognized in calculus textbooks that the fundamental theorem of calculus is simply a fixed-point equation:
\begin{equation}\label{FTOC2}
f = \LL f = f(a)\oo + \II \DD f,
\end{equation}
where \LL is defined in \eqref{ldef}.
According to Theorem~\ref{thm0}, \eqref{FTOC2} holds almost everywhere for all absolutely continuous functions $f: [a,b]\rightarrow \RR$.  

To obtain Taylor's theorem,  we apply a slight modification of the technique of fixed-point iteration described in Section~\ref{Fixed}.
If the derivative $\DD f$ is absolutely continuous, 
then by  \eqref{FTOC2}  we have $\DD f = \LL Df  = \DD f(a)\oo + \II \DD (\DD f)$. Making this replacement in  \eqref{FTOC2} gives
 \begin{equation}\label{iter2}
 \begin{aligned}
 f &=  f(a)\oo + \II \DD f\\
 &=  f(a)\oo + \II (\, \DD f(a)\oo + \II \DD (\DD f) \,)\\
 &= f(a)\oo +  \DD f(a)\II \oo + \II (\II \DD (\DD f))\\ 
 &= f(a)\oo + \DD f(a)\II  \oo + \II^2 \DD^2 f.
 \end{aligned}
 \end{equation}
 The last two equalities in \eqref{iter2} follow from the linearity and associativity properties of operators described in Section~\ref{ops}. 
 
We can now continue the iterative process. If $\DD^2 f$ is absolutely continuous, then by \eqref{FTOC2} we may similarly replace $\DD^2 f$ in \eqref{iter2} with  $ \DD^2 f(a)\oo + \II \DD (\DD^2 f)$.  
This gives
 \begin{equation}\label{iter3}
 \begin{aligned}
 f  
 &= f(a)\oo + \DD f(a)\II \oo + \II^2 (  \DD^2 f(a) + \II \DD (\DD^2 f) )\\ 
 &= f(a)\oo + \DD f(a)\II \oo + \DD^2 f(a)\II^2  \oo + \II^3 \DD^3 f.
 \end{aligned}
 \end{equation}
 Clearly we may continue the same process for $\DD^3 f, \DD^4 f, \ldots \DD^N f$ as long as all of these derivatives exist and are absolutely continuous (in fact, if $\DD^N f$ is absolutely continuous, then all of the lower-order derivatives will  be absolutely continuous as well). We may summarize the result as follows:
 \begin{equation}\label{TTN}
 f =  f(a)\oo + \DD f(a)\II \oo + \DD^2 f(a)\II^2  \oo + \ldots + \DD^N f(a)\II^N  1 + \II^{N+1} \DD^{N+1} f.
 \end{equation}
Equation \eqref{TTN} is actually Taylor's theorem in disguise.  To see this, we only need to rewrite the functions $\II^n \oo$ in conventional notation, and evaluate the integrals successively:
\begin{equation}\label{TTNb}
\begin{aligned}
 \II \oo(x) &=  \int_{a}^{x} dt_1 = (x-a);\\
 \II^2 \oo(x) &= \int_{a}^{x} \II \oo(t_1) dt_1 = \int_{a}^{x} (t_1-a) dt_1 = \frac{(x-a)^2}{2};\\
 \II^3 \oo(x) &= \int_{a}^{x} \II^2 \oo(t_1)  dt_1 = \int_{a}^{x} \frac{(t_1-a)^2}{2} dt_1 = \frac{(x-a)^3}{3!};\\
 \ldots &\ldots \ldots  \\
 \II^N \oo(x) &= \int_{a}^{x} \II^{N-1} \oo(t_1) dt_1 = \int_{a}^{x} \frac{(t_1-a)^{N-1}}{(N-1)!} dt_1 = \frac{(x-a)^N}{N!}.
 \end{aligned}
\end{equation}
It remains to evaluate the final term in \eqref{TTN}, which is $\II^{N+1} \DD^{N+1} f$. In this case, the  integrand is not the constant function $\oo$ as in the other integrals. However, we may use \eqref{Imon2} to give an upper bound on this term:
\begin{equation}
\begin{aligned}
| \II^{N+1} \DD^{N+1} f |  &\le \sup_{c\in [a,x]}| \DD^{N+1} f(c)| \II^{N+1}  \oo\\
& \le  \sup_{c\in [a,x]}| f^{(N+1)}(c)| \frac{(x-a)^{N+1}}{(N+1)!}
\end{aligned}
\end{equation}
which is the bound on $R_N$ in \eqref{boundR}. For the exact evaluation of $R_N$, 
 we may repeatedly use the fact that
integral order can be exchanged:
\begin{equation}\label{exch}	
\int_a^{t_k} \int_a^{t_j}  g(t_i,t_j)  dt_i\,dt_j = \int_a^{t_k}  \int_{t_i}^{t_k}  g(t_i,t_j)  dt_j\,dt_i.
\end{equation}
Notice what happens to the integral limits under exchange:  the limits of the outer integral do not change; the lower limit of the inner integral becomes the outer integration variable; and the upper limit of the inner integral is the same as the upper limit of the outer integral. We may apply this rule first
to exchange the $dt_N$ and $dt_{N+1}$ integrals (shown in parentheses):
\begin{equation}\label{TTNc0}
\begin{aligned}
  &\II^{N+1} \DD^{N+1} f \\
  &\quad = \int_{a}^{x}  \ldots \int_a^{t_{N-2}} \left( \int_a^{t_{N-1}} \int_a^{t_{N}}   f^{(n+1)}(t_{N+1})  dt_{N+1} dt_N \right)  dt_{N-1} \ldots  dt_1\\
   &\quad= \int_{a}^{x} \ldots \int_a^{t_{N-2}} \left(\int_a^{t_{N-1}}  \int_{t_{N+1}}^{t_{N-1}} f^{(n+1)}(t_{N+1})   dt_{N} dt_{N+1} \right)dt_{N-1}   \ldots   dt_1,
 \end{aligned} 
  \end{equation}   
  Next, we exchange the $dt_{N-1}$ and $dt_{N+1}$ integrals in similar fashion:
\begin{equation}\label{TTNc1}
\begin{aligned}
  &\II^{N+1} \DD^{N+1} f \\
   &\quad= \int_{a}^{x} \ldots\left( \int_a^{t_{N-2}} \int_a^{t_{N-1}}  \left[ \int_{t_{N+1}}^{t_{N-1}} f^{(n+1)}(t_{N+1})   dt_{N} \right]  dt_{N+1} dt_{N-1} \right)  \ldots   dt_1\\
   &\quad= \int_{a}^{x} \ldots\left( \int_a^{t_{N-2}} \int_{t_{N+1}}^{t_{N-2}}  \left[ \int_{t_{N+1}}^{t_{N-1}} f^{(N+1)}(t_{N+1})   dt_{N} \right]  dt_{N-1} dt_{N+1} \right)  \ldots   dt_1,
 \end{aligned} 
  \end{equation} 
where the integral in square brackets plays the role of $g(t_i,t_j)$ in \eqref{exch}, and the integrals in parentheses have been exchanged. In the same way we may exchange the $dt_{N+1}$ integral successively with $dt_{N-2},dt_{N-3}, \ldots dt_1$. The final result is 
\begin{equation}\label{TTNc}
\begin{aligned}  
 &\II^{N+1} \DD^{N+1} f \\
   &\quad= \int_{a}^{x}   \int_{t_{N+1}}^{t_{1}} \ldots \int_{t_{N+1}}^{t_{N-1}}   f^{(N+1)}(t_{N+1})   dt_{N} \ldots dt_{1}  dt_{N+1} \\
   &\quad= \int_{a}^{x}  f^{(N+1)}(t) \left(\int_t^{t_{1}} \ldots\int_{t}^{t_{N-1}}      dt_{N} \ldots dt_{1}\right)  dt, 
 \end{aligned}
\end{equation}
where the last  expression in \eqref{TTNc} is obtained by replacing $t_{N+1}$ with $t$, and by moving $f^{(n+1)}(t)$ outside of the integrals over  $dt_1,\ldots dt_N$. The $N$-fold integral in parentheses is identical to the final integral in \eqref{TTNb}, except that the lower limit is $t$ instead of $a$. It follows that
\begin{equation}\label{TTNd}
\begin{aligned}
  \II^{N+1} \DD^{N+1} f & =  \int_{a}^{x} f^{(n+1)}(t) \frac{(x-t)^N}{N!}  dt,
 \end{aligned}
\end{equation}
which agrees with \eqref{exactR}.

\subsection{Alternative evaluation using volume integrals}
An alternative, elegant evaluation of the integrals in \eqref{TTNc} is acheived by interpreting the integrals as volume integrals over regions in $\RR^n, 1\le n \le N$. 
The $n$-dimensional integral $\int_{a}^{x} \int_a^{t_1} \ldots \int_a^{t_n} dt_n \ldots dt_2 dt_1$ can be interpreted as the volume of the set $\AA_n$ in $\mathbb{R}^n$, where 
\begin{equation}
\AA_n := {a \le t_n \le t_{n-1} \le \ldots \le t_1 \le x }.
\end{equation}
But the ordering $t_n \le t_{n-1} \le \ldots \le t_1$ is simply one of $n!$ possible orderings of the $n$ variables $t_1,\ldots t_n$.  Since these variables are dummy variables that are integrated over, the value of the integral does not depend on the ordering of the variables.  For example, in the case where $n= 3$ we have $3!=6$ different orderings of the variables, namely:
\begin{equation}
\begin{aligned}
&a \le t_3 \le t_2 \le t_1 \le x; \quad a \le t_3 \le t_1 \le t_2 \le x; \quad  a \le t_2 \le t_3 \le t_1 \le x;\\
&a \le   t_2 \le t_1 \le t_3 \le x;\quad   a \le t_1 \le t_3 \le t_2 \le x; \quad  a \le t_1 \le t_2 \le t_3 \le x.
\end{aligned}
\end{equation} 
But these 3! orderings  correspond to 3! disjoint sets that together make up the cube $[a,x]
^3$ \footnote{There is a technical issue here in that the surfaces of these sets are not disjoint, but comprise a set of measure 0. However, it is intuitively clear that the volumes of these sets should add to the volume of the cube.  In three dimensions, this can be demonstrated using a model.}. The volumes of these $3!$ sets are equal:  therefore the volume of $\AA_3$  is $1/3!$ of the volume of the cube, giving the same result as \eqref{TTNb} for $N=3$.  
The argument generalizes directly to $n$ dimensions ($1 \le n \le N$), giving the result:
\begin{equation}\label{A_int}
\int_{a}^{x} \int_a^{t_1} \ldots \int_a^{t_n} dt_n \ldots dt_2 dt_1 = (\text{Volume of }\AA_n) = \frac{(x-a)^n}{n!},
\end{equation}
which is the same as \eqref{TTNb}. To evaluate the final integral \eqref{TTNd}, we note that
\begin{equation}
\AA_{N+1}  \cap \{t_{N+1} =t \} =  \AA'_{N}(t),
\end{equation}
where 
\begin{equation}
 \AA'_{N}(t) := {t \le t_N \le t_{N-1} \le \ldots \le t_1 \le x }
\end{equation}
so that the volume integral $\II^{N+1} \DD^{N+1} f $ can be expressed as:
\begin{equation}\label{lastInt}
\begin{aligned}
\II^{N+1} \DD^{N+1} f  &=  \int_a^x \left( \int_{\AA'_{N}(t)} f^{(N+1)}(t) dt_N \ldots dt_2 dt_1 \right)  dt\\
&=  \int_a^x f^{(N+1)}(t) \left( \int_{\AA'_{N}(t)} dt_N \ldots dt_2 dt_1\right)  dt\\
& = \int_a^x f^{(N+1)}(t) \frac{(t-a)^N}{N!}  dt,
\end{aligned}
\end{equation}
where we have used the volume formula \eqref{A_int} with $n \rightarrow N$ and $x \rightarrow t$. This result 
agrees with \eqref{TTNd}. 

\section{Discussion}\label{Discussion}
The role of mathematics within science and society is changing rapidly. Because of the digital revolution, math is having an increasingly weighty impact on all aspects of society, including business and government. This impact comes through the need for mathematically-based algorithms and procedures in communication, classification, modeling, prediction, and control. Many of these results are framed in terms of linear algebraic concepts such as vectors, matrices, and tensors, while calculus and functional analysis enter in cases where the systems under study can be approximated as continuous.
 
The idea of linear operator (or functional) as a function of functions is an important bridge between linear algebra and mathematical analysis. It is thus one of the cornerstones of modern applied mathematics. 
In order to adequately prepare students to deal with current mathematical challenges in science and technology, the teaching of mathematics should make appropriate adaptations.
 The proof of Taylor's theorem in Section~\ref{Method} emphasizes the relation between numerical methods and analytical theory, and is thus suitable for this purpose.
  
One possible objection to the proof is that it involves long equations (such as  \eqref{TTNc1}) that appear very complicated. 
Many students may find such long equations difficult to deal with, because they lack the technical facility in algebraic manipulations.  But in applied mathematics  (particularly areas of mathematics related to computation, such as linear algebra and numerical analysis) these types of equations are very common. Although they appear complicated, in actuality these equations are built up from the repeated application of very simple ideas.    Furthermore, such equations may readily handled by   looking at simple cases first, and then building up to the general case. The proof in Section~\ref{Method} is a good example of this approach. For this reason, we assert that the new perspective on Taylor's theorem provided in this paper can be a valuable addition to the undergraduate mathematics curriculum. 

\bibliographystyle{plain}
\bibliography{taylorBib}

\end{document}